\documentclass[a4paper, 12pt]{article}
\usepackage[utf8]{inputenc}
\usepackage{cmap}
\usepackage{amsfonts}
\usepackage[T2A]{fontenc}
\usepackage[english]{babel}
\usepackage{amsmath}
\usepackage{amssymb}
\usepackage{amsthm}
\usepackage[a4paper,hmargin=2.5cm,vmargin=2.5cm]{geometry}
\usepackage{euscript}
\usepackage{graphicx}
\usepackage{float}
\usepackage{wrapfig}
\usepackage{hyperref}

\DeclareMathOperator{\Rep}{Rep}

\newtheorem{theorem}{Theorem}
\newtheorem{prop}{Proposition}
\newtheorem{lemma}{Lemma}

\title{Stability of Restrictions of Representations of the Symmetric Group to the Hyperoctahedral Subgroup}

\date{}

\begin{document}

\maketitle

\begin{center}
    {\large Sergey Davydov}
\end{center}

\begin{abstract}
The paper investigates the stability properties of restrictions of irreducible representations of the symmetric group to the hyperoctahedral subgroup. A stability result is obtained, analogous to the classical Murnaghan theorem on the stability of the decomposition of tensor products of representations of the symmetric group. The proof is based on the description of these restrictions in terms of symmetric functions from \cite{K-T}.
\end{abstract}

\section{Introduction}

If $\nu = (\nu_1, \nu_2, \ldots)$ is a partition (Young diagram), then $\nu[n]$ denotes the partit.ion
$$\nu[n] = (n - |\nu|, \nu_1, \nu_2, \ldots).$$
In this notation, it is assumed that $n$ is sufficiently large so that $\nu[n]$ is a valid partition, meaning $n \ge |\nu| + \nu_1$.

Let $V^{\lambda}$ be the irreducible representation of the symmetric group $S_n$ (where $n = |\lambda|$) corresponding to the partition $\lambda$. A classical result concerning the stability of Kronecker coefficients states that for a pair of Young diagrams $\lambda, \mu$ there exists a decomposition
\begin{equation}
V^{\lambda[n]} \otimes V^{\mu[n]} = \bigoplus_{\nu} \gamma_{\lambda\mu}^\nu V^{\nu[n]},
\label{tens}
\end{equation}
which does not depend on $n$ for sufficiently large $n$.

This statement was formulated by Murnaghan in \cite{Murn}, but he did not provide a proof. A complete proof was later given by Littlewood in \cite{Littl} in the mid-20th century.

Since then, alternative proofs of Murnaghan's theorem and far-reaching generalizations of the fact have been proposed in numerous works. For example, in the paper~\cite{Thib1} Thibon uses the specifics of Hopf algebra structure and $\lambda$-structure on the ring of symmetric functions to derive certain stable formulas for the decomposition of tensor products~(\ref{tens}). In a more recent work~\cite{Thib2} the same author reformulates this approach in terms of certain vertex operators and obtains results on the stability of inner plethysm of representations (inner plethysm refers to the
$\lambda$-ring structure on the representation ring $\Rep(S_n)$, see \cite{Knut}).

A fully combinatorial proof of the stability of Kronecker coefficients is given by E.~Vallejo in~\cite{Val} where he expresses Kronecker coefficients in terms of Littlewood-Richardson coefficients and uses a combinatorial interpretation of the inverse Kostka matrix elements.

A completely different proof of stability is provided by S. Sam and A. Snowden in~\cite{Snow}. Their approach, using Schur-Weyl duality, relies on the Borel-Weil theorem, which establishes a correspondence between line bundles over flag varieties and representations of $GL_n$.

The problem of decomposing the tensor product of representations into irreducible components can be seen as a special case of the problem of decomposing the restriction of a representation: $V^{\lambda}\otimes V^{\mu}$ is the restriction of the representation $V^{\lambda} \boxtimes V^{\mu}$ (i.e., the tensor product considered as a representation of the group $S_n \times S_n$) to the diagonal subgroup $S_n\subset S_n\times S_n$.

From this perspective, a similar question can be posed in other contexts. For instance, the hyperoctahedral group
$$H_n = S_n \ltimes (\mathbb Z /2 \mathbb Z)^n$$
is naturally embedded in $S_{2n}$ as the centralizer of the permutation $(12)(34)\ldots(2n-1,2n)$. This allows us to consider the restriction of irreducible representations of the symmetric group $S_{2n}$ to $H_n$ and raises the natural question: does a stability phenomenon of the same kind occur in this case, and if so, in what form?

As it turns out, such stability indeed exists, and this is demonstrated in the present work.

\section{Formulation and Preliminaries}

The irreducible representations of the hyperoctahedral group $H_n$
are parameterized by pairs of partitions $(\lambda, \mu)$ such that $|\lambda| + |\mu| = n$ (where $\lambda$ corresponds to the trivial representation of $\mathbb Z /2 \mathbb Z$ and $\mu$ corresponds to the sign representation; for details on how the pairs parameterize representations, see~\cite{Serre}, $\S9.2$). These representations are denoted by $V^{\lambda,\mu}$. Let $l(\lambda)$ denote the number of nonzero parts in $\lambda$, and let $\lambda'$ denote the conjugate (transposed) partition of $\lambda$.

\begin{theorem}
Let $\nu$ be a partition. Then, for sufficiently large $n$, there exists a decomposition
$$
Res_{H_n}^{S_{2n}} \left(V^{\nu[2n]}\right) = \bigoplus_{\lambda,
\mu} \overline{K}_{\lambda, \mu}^{\nu}V^{\lambda[n-|\mu|],\ \mu},
$$
where the coefficients $\overline{K}_{\lambda, \mu}^{\nu}$ do not depend on $n$, and the summation is taken over a certain finite set of pairs of partitions $\lambda, \mu$.
\label{main}
\end{theorem}

The main tool for proving this theorem is the following fact, established in \cite{K-T} (Proposition 2.2).

\begin{prop}
The multiplicity $K_{\lambda,\mu}^{\nu}$ of the representation $V^{\lambda, \mu}$ in the restriction $Res_{H_n}^{S_{2n}}(V^{\nu})$ is equal to the coefficient of $s_{\nu}$ in the expansion of the function $$(s_{\lambda}\circ h_2)(s_{\mu}\circ e_2)$$
in terms of the basis of Schur functions, where $"\circ"$ denotes the (outer) plethysm of symmetric functions, and $h_2$ and $e_2$ are the complete and elementary symmetric functions of degree $2$ respectively.
\label{reduction}
\end{prop}

In addition to standard information about symmetric functions (such as the Jacobi-Trudi formula and the Littlewood-Richardson rule), we will also use the following formulas.

\begin{prop}
For plethysm, the following identities hold:
\end{prop}
    \begin{equation}
        h_k\circ h_2 = \sum_{|\alpha| = k} s_{2\alpha}, \label{a}
    \end{equation}
{\itshape where  $2\alpha = (2\alpha_1, 2\alpha_2, \ldots).$}
    \begin{equation}
        e_k\circ h_2 = \sum_{\substack{|\beta| = k \\ \beta_1 > \beta_2 > \ldots} }s_{\Gamma(\beta)'}, \label{b}
    \end{equation}
{\itshape where $\Gamma(\beta) = (\beta_1 - 1, \beta_2 - 1, \ldots| \beta_1, \beta_2, \ldots)$ in Frobenius notation,}
    \begin{equation}
        h_k\circ e_2 = \sum_{|\alpha| = k} s_{(2\alpha)'}.\label{c}
    \end{equation}
    \begin{equation}
        e_k\circ e_2 = \sum_{\substack{|\beta| = k \\ \beta_1 > \beta_2 > \ldots} }s_{\Gamma(\beta)}.\label{d}
    \end{equation}

Proofs of these formulas can be found in the monograph~\cite{M} (Chapter 1, Section 8, Example 6). A more detailed proof of formulas~(\ref{b}) and~(\ref{d}) can also be found in \cite{K-T} (Lemma 1.1).

\section{Main Lemmas}
Let $\Lambda$ denote the ring of symmetric functions.

We introduce a useful function on $\Lambda$:
$$
\Phi \left(\sum_{\lambda}c_{\lambda} s_{\lambda}\right) = \max_{\lambda:\ c_{\lambda} \neq 0} \lambda_1
$$
It is obvious that $\Phi(f) \le \deg f$. Let us also note the following important property of the function $\Phi$.

\begin{prop}
$\Phi(fg) \le \Phi(f) + \Phi(g)$.
\end{prop}
\textbf{Proof.}
It is sufficient to prove this for the case $f = s_{\lambda},\ g = s_{\mu}$. In this case, the inequality follows from the Littlewood-Richardson rule, since if $\nu$ is a partition such that $\nu_1 > \lambda_1 + \mu_1$, then the Littlewood-Richardson coefficient $c_{\lambda \mu}^{\nu}$ is zero. This follows from the fact that in the skew diagram $\nu/\lambda$ the first row contains more than $\mu_1$ cells, meaning that it is impossible to fit a reverse lattice word with $\mu_1$ ones. $\square$

\begin{lemma}
    Fix three partitions $\lambda, \mu, \nu$. Then the coefficient $K_{\lambda, \mu}^{\nu}(n) := K_{\lambda[n-|\mu|], \mu}^{\nu[2n]}$ does not depend on $n$ for $2n \ge |\nu| + \nu_1 + 2|\lambda| + 2|\mu|$.
    \label{coef_stability}
\end{lemma}

\textbf{Proof.} According to Proposition~\ref{reduction}, 
\begin{equation}
K_{\lambda, \mu}^{\nu}(n) = \left\langle(s_{\lambda[n-|\mu|]}\circ h_2)(s_{\mu}\circ e_2), s_{\nu[2n]}\right\rangle,
\label{scal}
\end{equation}
where $\langle\cdot , \cdot\rangle$ denotes the standard inner product on the ring of symmetric functions (Schur functions form an orthonormal basis). We express $s_{\lambda[n-|\mu|]}$ using the Jacobi-Trudi formula: $s_{\gamma} = \det (h_{\gamma_i + j - i})$, expanding the determinant along the first row:
\begin{equation}
s_{\lambda[n - |\mu|]} = \sum_{i=0}^{l(\lambda)} h_{n - |\mu| - |\lambda| + i}\cdot f_i,
\label{Jac-Tr}
\end{equation}
where $f_i$ are certain polynomials with $\deg f_i = |\lambda| - i \le |\lambda|$.

Next, applying plethysm $\circ h_2$ to~(\ref{Jac-Tr}), using the fact that plethysm is a homomorphism with respect to the left argument, we obtain:

\begin{equation}
s_{\lambda[n - |\mu|]}\circ h_2 = \sum_{i=0}^{l(\lambda)} (h_{n - |\mu| - |\lambda| + i}\circ h_2)\cdot g_i,
\label{Jac-Tr2}
\end{equation}
where $g_i = f_i\circ h_2$, and $\deg g_i = 2\deg f_i \le 2|\lambda|.$ Now, multiplying by the polynomial $s_{\mu}\circ e_2$ of degree $2|\mu|$ we get:
\begin{equation}
(s_{\lambda[n-|\mu|]}\circ h_2)(s_{\mu}\circ e_2) = \sum_{i=0}^{l(\lambda)} (h_{n - |\mu| - |\lambda| + i}\circ h_2)\cdot r_i,
\label{Jac-Tr3}
\end{equation}
where $\deg r_i \le 2|\lambda| + 2|\mu|$.

Next, we expand the plethysms $h_k\circ h_2$ in~(\ref{Jac-Tr3}) using equation~(\ref{a}), yielding Schur polynomials
$s_{2\alpha}$ for $\alpha \vdash n - |\mu| - |\lambda| + i$,
which are multiplied by polynomials $r_i$ of degree at most $2|\lambda| + 2|\mu|$. Since in~(\ref{scal}) we take the inner product with $s_{\nu[2n]}$, and
$$\Phi(s_{\nu[2n]}) = 2n - |\nu|,$$
only those terms $s_{2\alpha} \cdot r_i$ remain where
$$\Phi(s_{2\alpha} \cdot r_i) \ge 2n - |\nu|.$$
However
$$\Phi(s_{2\alpha}r_i) \le \Phi(s_{2\alpha}) + \Phi(r_i) \le 2\alpha_1 + 2|\lambda| + 2|\mu|.$$
Thus, we only need to consider terms $s_{2\alpha}r_i$ where
$$2\alpha_1 + 2|\lambda| + 2|\mu| \ge 2n - |\nu|.$$
Furthermore, since
$$
|\alpha| = n - |\mu| - |\lambda| + i \le n - |\mu| - |\lambda| + l(\lambda),
$$
it follows that 
$$|2\alpha| - 2\alpha_1 \le (2n - 2|\lambda| - 2|\mu| + 2l(\lambda)) - (2n - 2|\lambda| - 2|\mu| - |\nu|) = 2l(\lambda) + |\nu|,$$
That is, the number and type of such partitions do not depend on $n$. Expanding each $r_i$ in the Schur function basis $s_{\beta}$, we need to analyze the inner product
\begin{equation}
\label{scal2}
\left\langle \sum_{\kappa,\beta} c_{\kappa\beta}s_{\kappa}s_{\beta}, s_{\nu[2n]}\right\rangle,
\end{equation}
where $\kappa_1 \ge 2n - 2|\lambda| - 2|\mu| - |\nu|$ (with $\kappa = 2\alpha$). According to the Littlewood-Richardson rule, we only need to consider those $\kappa$ for which $\kappa \backslash \kappa_1 \subset \nu$. As noted earlier, in formula~(\ref{scal2}) the summation runs over a fixed set of pairs of partitions $\kappa, \beta$ with fixed coefficients $c_{\kappa\beta}$, and the only part of $\kappa$ that depends on $n$ is its first row, which consists of
$$2n-|\beta| - \sum_{i=2}\kappa_i$$
cells. The Littlewood–Richardson coefficient $c_{\kappa\beta}^{\nu[2n]}$ does not depend on $n$ if $\kappa_1 \ge \nu_1$, because the skew diagram $\nu[2n]/\kappa$ does not depend on $n$, up to a shift of its connected components.

For this inequality to hold, it is sufficient that
$$2n - 2|\lambda| - 2|\mu| - |\nu| \ge \nu_1,$$
This is exactly the condition stated in the lemma. $\square$

From Lemma~\ref{coef_stability} the main theorem does not yet follow because the estimate in the lemma depends on the diagrams $\lambda$ and $\mu$. Now we will show that $K^{\nu}_{\lambda,\mu}(n) = 0\ \forall n$ for all pairs $\lambda, \mu$, except for a finite number.

\begin{lemma}
    If $K^{\nu[2n]}_{\lambda, \mu} \neq 0$, then $|\lambda| - \lambda_1 \le |\nu|$ and $|\mu| \le |\nu|.$
    \label{finiteness}
\end{lemma}
\textbf{Proof.} Suppose that
$$
\left\langle (s_{\lambda}\circ h_2)(s_{\mu}\circ e_2), s_{\nu[2n]}\right\rangle \neq 0.
$$
Expanding $s_{\lambda}$ using the dual Jacobi-Trudi formula: $s_{\lambda} = \det(e_{\lambda_j' + i - j})$, we obtain
$$s_{\lambda}\circ h_2 = \det(e_{\lambda_j' + i - j}\circ h_2),$$
This determinant  is an alternating sum of products of the form $$\prod_{j=1}^{l(\lambda') = \lambda_1} e_{k_j}\circ h_2$$
with $\sum_j k_j = |\lambda|$, and for each such product we have the estimate
$$
\Phi\left(\prod_j^{\lambda_1} e_{k_j}\circ h_2\right) \le \sum_{j=1}^{l(\lambda') = \lambda_1}\Phi(e_{k_j}\circ h_2) = \sum_{j=1}^{\lambda_1} (k_j +1) = |\lambda|+\lambda_1.
$$
Here, we used the fact that
$$\Phi(e_{k}\circ h_2) = k+1,$$
which follows from formula~(\ref{b}) (indeed, $\Phi(s_{\Gamma(\beta)'}) = \beta_1+1$). At the same time,
$$\Phi(s_{\mu}\circ e_2) \le \deg (s_\mu\circ e_2) = 2|\mu|,$$
which implies
$$\Phi((s_{\lambda}\circ h_2)(s_{\mu}\circ e_2)) \le |\lambda|+\lambda_1+2|\mu| = 2n - |\lambda| +\lambda_1$$
(since $|\lambda| + |\mu| = n$). On the other hand, from the condition of the lemma, $$\Phi((s_{\lambda}\circ h_2)(s_{\mu}\circ e_2)) \ge 2n - |\nu|.$$
Thus, we obtain
$$2n - |\lambda| +\lambda_1 \ge 2n - |\nu|$$
which simplifies to $|\lambda| - \lambda_1 \le |\nu|$.

A similar argument establishes the second inequality. Expanding $s_{\mu}$ using the Jacobi-Trudi formula: $s_{\mu} = \det(h_{\mu_i+j-i})$, we obtain 
$$\Phi(s_{\mu}\circ e_2) = \Phi(\det (h_{\mu_i+j-i}\circ e_2)) \le \max_{k_1 + \ldots + k_{l(\mu)} = |\mu|}\left(\sum_{j=1}^{l(\mu)}\Phi(h_{k_j}\circ e_2)\right) = \sum_{j=1}^{l(\mu)} k_j =  |\mu|,$$
where we used the fact that
$$\Phi(h_k\circ e_2) = k,$$ which follows from formula~(\ref{c}) (since $\Phi(s_{(2\alpha)'}) = \alpha_1'$). Thus, we obtain
$$\Phi((s_{\lambda}\circ h_2)(s_{\mu}\circ e_2)) \le 2|\lambda| + |\mu| = 2n - |\mu|$$
On the other hand, we have
$$\Phi((s_{\lambda}\circ h_2)(s_{\mu}\circ e_2)) \ge 2n - |\nu|.$$
From the last two inequalities, we obtain $2n - |\mu| \ge 2n - |\nu|$, which simplifies to $|\mu| \le |\nu|$.~$\square$

\section{Conclusion and the case of diagrams with a long first column.}

From Lemmas~\ref{coef_stability} and~\ref{finiteness}, Theorem~\ref{main} follows directly. Moreover, by combining the bounds given in the lemmas, we obtain that the decomposition in Theorem~$\ref{main}$ necessarily stabilizes when
$$2n \ge 5|\nu| + \nu_1.$$
Additionally, by applying the isometric involution
$$\omega: \Lambda \to \Lambda,\ \omega(s_{\lambda}) = s_{\lambda'}$$
to formula~(\ref{scal}) and using the fact that
$$
\omega(f\circ h_2) = f\circ e_2 \ \forall f\in \Lambda,
$$
(which holds for $f = h_k$ due to formulas (\ref{a}) и (\ref{c}), and for an arbitrary $f$ follows from the fact that plethysm is a homomorphism with respect to the left argument, and $h_k$'s generate $\Lambda$) we obtain an analogous result for the conjugate diagrams $\nu[2n]':$
\begin{theorem}
Let $\nu$ be a partition. Then there exists a stable decomposition
$$
Res_{H_n}^{S_{2n}} \left(V^{\nu[2n]'}\right) = \bigoplus_{\lambda, \mu} \overline{K}_{\lambda, \mu}^{\nu}V^{\mu,\lambda[n-|\mu|]}.
$$
\end{theorem}

\section{Acknowledgements}

I would like to express my deep gratitude to my supervisor, Prof. Grigori Olshanski, for his guidance throughout this research. I also thank the staff at the Faculty of Mathematics, HSE, and the Scientific School named after I.M. Krichever, Skolkovo Institute of Science and Technology.

\end{document}